%
\catcode`@=11
%
%
\def\bibn@me{R\'ef\'erences}
\def\bibliographym@rk{\centerline{{\sc\bibn@me}}
    \sectionmark\section{\ignorespaces}{\unskip\bibn@me}
    \bigbreak\bgroup
    \ifx\ninepoint\undefined\relax\else\ninepoint\fi}
%
%
%
\let\refsp@ce=\
\let\bibleftm@rk=[
\let\bibrightm@rk=]
%
%
%
\def\numero{n\raise.82ex\hbox{$\fam0\scriptscriptstyle o$}~\ignorespaces}
%
%
\newcount\equationc@unt
\newcount\bibc@unt
\newif\ifref@changes\ref@changesfalse
\newif\ifpageref@changes\ref@changesfalse
\newif\ifbib@changes\bib@changesfalse
\newif\ifref@undefined\ref@undefinedfalse
\newif\ifpageref@undefined\ref@undefinedfalse
\newif\ifbib@undefined\bib@undefinedfalse
\newwrite\@auxout
%
%
\def\eqnum{\global\advance\equationc@unt by 1%
\edef\lastref{\number\equationc@unt}%
\eqno{(\lastref)}}
%
%
%
%
%
%
\def\re@dreferences#1#2{{%
    \re@dreferenceslist{#1}#2,\undefined\@@}}
\def\re@dreferenceslist#1#2,#3\@@{\def\next{#2}%
    \expandafter\ifx\csname#1@@\meaning\next\endcsname\relax
    ??\immediate\write16
    {Warning, #1-reference "\next" on page \the\pageno\space
    is undefined.}%
    \global\csname#1@undefinedtrue\endcsname
    \else\csname#1@@\meaning\next\endcsname\fi
    \ifx#3\undefined\relax
    \else,\refsp@ce\re@dreferenceslist{#1}#3\@@\fi}
%
%
%
\def\newlabel#1#2{{\def\next{#1}\newl@bel#2}}
\def\newl@bel#1#2{%
    \expandafter\xdef\csname ref@@\meaning\next\endcsname{#1}%
    \expandafter\xdef\csname pageref@@\meaning\next\endcsname{#2}}
\def\label#1{{%
    \toks0={#1}\message{ref(\lastref) \the\toks0,}%
    \ignorespaces\immediate\write\@auxout%
    {\noexpand\newlabel{\the\toks0}{{\lastref}{\the\pageno}}}%
    \def\next{#1}%
    \expandafter\ifx\csname ref@@\meaning\next\endcsname\lastref%
    \else\global\ref@changestrue\fi%
    \newlabel{#1}{{\lastref}{\the\pageno}}}}
\def\ref#1{\re@dreferences{ref}{#1}}
\def\pageref#1{\re@dreferences{pageref}{#1}}
%
%
\def\bibcite#1#2{{\def\next{#1}%
    \expandafter\xdef\csname bib@@\meaning\next\endcsname{#2}}}
\def\cite#1{\bibleftm@rk\re@dreferences{bib}{#1}\bibrightm@rk}
%
%
\def\beginthebibliography#1{\bibliographym@rk
    \setbox0\hbox{\bibleftm@rk#1\bibrightm@rk\enspace}
    \parindent=\wd0
    \global\bibc@unt=0
    \def\bibitem##1{\global\advance\bibc@unt by 1
        \edef\lastref{\number\bibc@unt}
        {\toks0={##1}
        \message{bib[\lastref] \the\toks0,}%
        \immediate\write\@auxout
        {\noexpand\bibcite{\the\toks0}{\lastref}}}
        \def\next{##1}%
        \expandafter\ifx
        \csname bib@@\meaning\next\endcsname\lastref
        \else\global\bib@changestrue\fi%
        \bibcite{##1}{\lastref}
        \medbreak
        \item{\hfill\bibleftm@rk\lastref\bibrightm@rk}%
        }
    }
\def\endthebibliography{\egroup\par}
%
%
%
\def\@closeaux{\closeout\@auxout
    \ifref@changes\immediate\write16
    {Warning, changes in references.}\fi
    \ifpageref@changes\immediate\write16
    {Warning, changes in page references.}\fi
    \ifbib@changes\immediate\write16
    {Warning, changes in bibliography.}\fi
    \ifref@undefined\immediate\write16
    {Warning, references undefined.}\fi
    \ifpageref@undefined\immediate\write16
    {Warning, page references undefined.}\fi
    \ifbib@undefined\immediate\write16
    {Warning, citations undefined.}\fi}
%
%
\immediate\openin\@auxout=\jobname.aux
\ifeof\@auxout \immediate\write16
  {Creating file \jobname.aux}
\immediate\closein\@auxout
\immediate\openout\@auxout=\jobname.aux
\immediate\write\@auxout {\relax}%
\immediate\closeout\@auxout
\else\immediate\closein\@auxout\fi
%
%
\input\jobname.aux
\immediate\openout\@auxout=\jobname.aux
%
%

\def\bibn@me{R\'ef\'erences bibliographiques}
%
\def\bibliographym@rk{\bgroup}
%
%
\outer\def\bye{     \par\vfill\supereject\end}

 \def\Leb{{\rm Leb}}
\overfullrule=0pt

\magnification=1200

  \def\pro{\noindent {\bf{Proof : }}}

\def\house#1{\setbox1=\hbox{$\,#1\,$}%
\dimen1=\ht1 \advance\dimen1 by 2pt \dimen2=\dp1 \advance\dimen2 by 2pt
\setbox1=\hbox{\vrule height\dimen1 depth\dimen2\box1\vrule}%
\setbox1=\vbox{\hrule\box1}%
\advance\dimen1 by .4pt \ht1=\dimen1
\advance\dimen2 by .4pt \dp1=\dimen2 \box1\relax}

  \def\eps{{\varepsilon}}

  \def\noi{\noindent}

\def\build#1_#2^#3{\mathrel{\mathop{\kern 0pt#1}\limits_{#2}^{#3}}}

\def\date {le\ {\the\day}\ \ifcase\month\or
janvier\or fevrier\or mars\or avril\or mai\or juin\or juillet\or
ao\^ut\or septembre\or octobre\or novembre\or
d\'ecembre\fi\ {\oldstyle\the\year}}

\font\fivegoth=eufm5 \font\sevengoth=eufm7 \font\tengoth=eufm10

\newfam\gothfam \scriptscriptfont\gothfam=\fivegoth
\textfont\gothfam=\tengoth \scriptfont\gothfam=\sevengoth

\def\cqfd{\unskip\kern 6pt\penalty 500 \raise 0pt\hbox{\vrule\vbox
to6pt{\hrule width 6pt \vfill\hrule}\vrule}\par}

\def\pro{\noindent {\it Proof. }}

\def\smallsquare{\vbox{\hrule\hbox{\vrule height 1 ex\kern 1 ex\vrule}\hrule}}
\def\cqfd{\hfill \smallsquare\vskip 3mm}

\def\Leb{{\lambda}}


\vskip 5mm

\centerline{\bf On fractional parts of powers
of real numbers close to $1$}

\vskip 13mm

\centerline{Yann B{\sevenrm UGEAUD} and Nikolay M{\sevenrm
OSHCHEVITIN} \footnote{}{\rm 2000 {\it Mathematics Subject
Classification: 11K31}.
 The second author was supported by
the RFBR grant No. 01-09-00371a }}

{\narrower\narrower
\vskip 15mm

\proclaim Abstract. {We prove that there exist arbitrarily small
positive real numbers $\eps$ such that every integral power
$(1 + \eps)^n$ is at a distance greater than $2^{-17}
\eps |\log \eps|^{-1}$ to the set of rational integers.
This is sharp up to the factor $2^{-17} |\log \eps|^{-1}$.
We also establish that the set of real numbers $\alpha > 1$
such that the sequence of fractional parts
$(\{\alpha^n\})_{n \ge 1}$ is not dense modulo $1$
has full Hausdorff dimension.}

}

\vskip 6mm

\vskip 12mm

\centerline{\bf 1. Introduction}

\vskip 6mm

Throughout this note, $\{ \cdot \}$ denotes the fractional
part and $|| \cdot ||$ the distance to the nearest integer.
In 1935, Koksma \cite{Ko35} established
that the sequence $(\{\alpha^n\})_{n \ge 1}$
is uniformly distributed modulo $1$,
for almost all
(with respect to the Lebesgue measure)
real numbers $\alpha$ greater than $1$.
However, very little is known on the distribution
of $(\{\alpha^n\})_{n \ge 1}$ for a specific real number $\alpha$
greater than $1$. If $\alpha$ is a Pisot number, that is,
an algebraic integer greater than $1$
all of whose Galois conjugates except
$\alpha$ lie in the open unit disc, then $||\alpha^n||$
tends to $0$ as $n$ tends to infinity, and the
limit points of $(\{\alpha^n\})_{n \ge 1}$
are contained in $\{0, 1\}$.
Pisot and Salem \cite{PiSa64} established that
if $\alpha$ is a Salem number, that is,
an algebraic integer greater than $1$
all of whose Galois conjugates except
$\alpha$ and $1/\alpha$ lie on the unit circle, then
$(\{\alpha^n\})_{n \ge 1}$ is dense but not uniformly
distributed modulo $1$. We do not know any explicit
transcendental real number $\alpha$ larger than $1$ for which
the sequence $(\{\alpha^n\})_{n \ge 1}$ is not uniformly
distributed modulo $1$.

In the present note, we are concerned with the set
$E$ composed of the real numbers $\alpha > 1$ for which
$(\{\alpha^n\})_{n \ge 1}$ is not dense modulo $1$.
In 1948, Vijayaraghavan \cite{Vij48} established that, for every
real numbers $a$ and $b$ with $1 < a < b$, the intersection
$E \cap (a, b)$ is uncountable.
Noticing that, in the proof of his Theorem 2, the parameter
$\eta$ should be taken equal to $\delta/(1 + b + \ldots + b^{h-1})$
and not to $1/(1 + b + \ldots + b^{h-1})$, the
following quantitative statement follows from his proof.

\proclaim Theorem V1.
There exist
arbitrarily small positive real numbers $\eps$ such that
$$
\inf_{n \ge 1} \, \Vert (1 + \eps)^n \Vert
>   \eps^{2/ \eps}.
$$

In the same paper,
Vijayaraghavan \cite{Vij48} also showed that, for any
interval $I$ of positive length contained in $[0, 1]$,
there are uncountably many $\alpha$ all of whose integral
powers are lying in $I$ modulo $1$. This result was recently
reproved by Dubickas \cite{Dub07c}.
Theorem 1 of \cite{Vij48}
includes the following statement.

\proclaim Theorem V2.
Let $H \ge 3$ be an integer. For every $\delta > 2/H$
and every interval $I$ of length $\delta$, there
exists $\alpha$ in $(H, H+1)$ such that $\{\alpha^n\}$
lies in $I$ for every $n \ge 1$.

The first purpose of the present note is to
significantly improve Theorem V1,
by means of a suitable modification of
a method introduced by Peres and Schlag \cite{PeSc09}
(see also \cite{Mosh09a,Mosh09b}),
based on the Lov\'asz local lemma.

Surprisingly, it seems that no metric result is known on the
size of the set $E$. The second aim of this note
is to give a suitable adaptation
of Vijayaraghavan's proof of
Theorem V2 for showing that $E$ has full Hausdorff dimension.

\vskip 5mm

\centerline{\bf 2. Main results}

\vskip 6mm

Our first result is a considerable improvement of Theorem V1.

\proclaim Theorem 1.
There exist arbitrarily small positive real numbers $\eps$
such that
$$
\inf_{n \ge 1} \, \Vert (1 + \eps)^n \Vert
> 2^{-17}  \, \eps  \, |\log \eps|^{-1}.
$$

Theorem 1 is sharp up to the
factor $2^{-17} |\log \eps|^{-1}$, since
the above infimum is clearly
at most equal to $\eps$, when $\eps < 1/2$.
The numerical constant $2^{-17}$ occurring in Theorem 1
can certainly be reduced, but we have made no effort in this direction.

The Peres--Schlag method is an inductive construction.
Roughly speaking, at each step $k$, we remove finitely many intervals,
which have (in all known applications until now)
essentially the same length. The novelty in the present application
of the method is that these intervals are far from having
the same length: here, at step $k$,
the quotient of the longest length by
the smallest one grows exponentially in $k$.
Consequently, the original approach of Peres and Schlag does
not allow us to prove Theorem 1, and we have to perform
a more complicated induction.

In \cite{BuMo11} we have combined the Peres--Schlag method
with the mass distribution principle to show that, in many
situations, the exceptional set constructed by means of the
Peres--Schlag method has full Hausdorff dimension.
A similar approach allows us to establish that,
for every small positive $\eps$,
the Hausdorff dimension of
$$
\bigl\{\eps' \in (\eps, 2 \eps) :
\inf_{n \ge 1} \, \Vert (1 + \eps')^n \Vert
> c  \, \eps  \, |\log \eps|^{-1} \bigr\}
$$
tends to $1$ when $c$ tends to $0$.
Brief explanations are given at the end of the proof
of Theorem 1.

The proof of Theorem 1 can be readily adapted to give the
more general following statement.

\proclaim Theorem 2.
Let $M$ be a positive real number.
For any non-zero real number $\xi$ in $[-M, M]$ and for any
sequence $(\eta_n)_{n \ge 1}$ of real numbers,
there exist a positive number $\gamma$,
depending only on $M$,
and arbitrarily small positive real numbers $\eps$ such that
$$
\inf_{n \ge 1} \, \Vert \xi (1 + \eps)^n + \eta_n \Vert
> \gamma  \, \eps  \, |\log \eps|^{-1}.
$$

Our last result
implies that the set of real numbers greater than 1
all of their integral powers stay, modulo one, in a given interval
of positive length
is rather big. It strengthens Corollary 5 of \cite{Dub07c}.

\proclaim Theorem 3.
Let $\xi$ be a positive real number.
Let $\eps < 1$ be a positive real number.
Let $(a_n)_{n \ge 1}$ be a sequence of real numbers
satisfying $0 \le a_n < 1 - \eps$ for $n \ge 1$.
The set of real numbers $\alpha$
such that $a_n \le \{\xi \alpha^n\} \le a_n + \eps$
for every $n \ge 1$ has
full Hausdorff dimension.

Theorem V2 suggests the next question,
which seems to be quite difficult.

\proclaim Question.
Let $\eps$ be a positive real number.
Are there arbitrarily large real numbers
$\alpha$ such that $\alpha$
is not a Pisot number
and all the fractional parts $\{\alpha^n\}$,
$n \ge 1$, are lying in an interval
of length $\eps / \alpha$ ?

Dubickas \cite{Dub07c} gave an alternative proof
of a version of Theorem V2 in which the lower
bound $2/H$ is replaced by $8/H$.

Throughout the present paper,
$\Leb$ denotes the Lebesgue measure.
Furthermore, $\lfloor x \rfloor$ and $\lceil x \rceil$
denote respectively the largest integer smaller than
or equal to $x$ and the smallest integer greater than or
equal to $x$.

\vskip 5mm

\centerline{\bf 3. Proof of Theorem 1}

\vskip 6mm

First, note that if the real numbers $\eps, \delta$
and the positive integers $k, m$
satisfy $0 < \eps, \delta < 1/5$ and
$$
|(1 + \eps)^k - m | \le \delta,
$$
then we get
$$
\biggl| {(1+\eps)^k \over m} - 1 \biggr| \le {\delta \over m}
$$
and
$$
\biggl| \log (1 + \eps) - {\log m \over k}  \biggr|
\le {2 \delta \over k m}.
$$

\bigskip

\noi{\it 1. Dangerous sets.}

\medskip

Let $t$ be a large positive integer and set
$$
\eta = 2^{-t}, \qquad \psi =  {\eta \over 2^{14} \log (1/\eta )}
= {1 \over 2^{t+14} t \log 2}. \eqno (3.1)
$$
Our aim is to find a real number $\xi$ such that
$$
2^{-t} \le \xi \le 2^{-t+1}
$$
and, for every $k \ge 1$,
$$
\xi  \not\in   A_k =
\bigcup_{m=\lfloor e^{\eta k} \rfloor}^{\lceil e^{2 \eta k} \rceil}
A_{k,m}, \quad
\hbox{where} \quad A_{k,m}
=  \left( { \log m \over k} -  {\psi \over k m} ,
{\log m \over k} + {\psi \over k m} \right).  \eqno (3.2)
$$
Setting $\eps := e^{\xi} - 1$, this proves our theorem
in view of the preliminary observation. Indeed, $\eps$
then satisfies
$$
|(1+\varepsilon)^k - m| \ge {\psi \over 2},
$$
for every positive integers $k, m$.

In the union occurring in (3.2), the integer $m$ varies between
$\lfloor e^{\eta k} \rfloor$ and
$\lceil e^{2 \eta k} \rceil$. Since the quotient
of these two numbers depends on $k$, we cannot use the
Peres--Schlag method as it was
applied in \cite{PeSc09,Mosh09a,Mosh09b}.
Fortunately, it
is possible to adapt it to prove our theorem.

In the sequel, we use an inductive process to establish that
$$
[2^{-t}, 2^{-t+1}] \setminus
\bigcup_{k \ge 1}
\bigcup_{m=\lfloor e^{\eta k} \rfloor}^{\lceil e^{2 \eta k} \rceil}
A_{k,m}
$$
is non-empty.
Put
$$
h = (t+6) 2^{t+6} = {2^6 \over \eta \log 2} \, \log {2^6 \over \eta}
\eqno (3.3)
$$
and
$$
k_n = h n, \quad \hbox{for $n \ge 1$}.
$$

\bigskip

\noi{\it 2. Initial steps.}

\medskip

We construct real numbers
$$
W_1,W_2, \ldots , W_{n-1}, \ldots  \in[\eta, 2\eta],
\qquad W_n  \le W_{n+1}, \qquad n \ge 1,
$$
and positive integers
$$
l_1,l_2, \ldots ,l_{n}, \ldots ,
\quad
w_1, w_2, \ldots , w_n, \ldots
$$
in such a way that
$$
W_n =  {w_n \over 2^{l_n}}, \qquad n \ge 1,
$$
and
$$
J_n =\left[W_n , W_{n} + {1 \over 2^{l_{n}}}\right] \subset
J_{n-1} =\left[W_{n-1} , W_{n-1} + {1 \over 2^{l_{n-1}}}\right]
\subset [\eta, 2\eta], \qquad n \ge 2.
$$

Let $l_1$ be such that
$$
2^{-l_1} < {2 \psi \over h \lceil e^{2 \eta h} \rceil}
\le 2^{-l_1 + 1},
$$
and observe that
$$
l_1 \ge 5 t,  \eqno (3.4)
$$
if $t$ is large enough.
Then, for each set $A_{k,m}$ with $k \le k_1$,
we consider the shortest dyadic interval $\hat{A}_{k,m}^1$
of the form
$$
\left( {a_1  \over 2^{l_1}},  {a_2 \over 2^{l_1}}
\right), \qquad a_1, a_2 \in \bf{Z},
$$
which covers the interval $A_{k,m}$, and we
define
$$
\hat{A}_k^1 :=
\bigcup_{m=\lfloor e^{\eta k} \rfloor}^{\lceil e^{2 \eta k} \rceil}
\hat{A}_{k,m}^1 \supset A_k.
$$
The choice of $l_1$ implies that
$$
\Leb (\hat{A}_{k,m}^1) \le 4 \Leb ({A}_{k,m}),
$$
for $\lfloor e^{\eta k} \rfloor \le m \le \lceil e^{2 \eta k} \rceil$.

Furthermore,
$$
\sum_{k=1}^h \,
\sum_{m = \lfloor e^{ \eta k}
\rfloor}^{\lceil e^{2 \eta k} \rceil} \, {2 \psi \over mk}
\le \psi h 2^{-t+3} \le 2^{-5} 2^{-t}.
$$
By (3.4),
this shows that there exists
$J_1 := [W_1, W_1 + 2^{-l_1}]$ such that
$ J_1 \cap \hat{A}_{k}^1 = \emptyset $
for every $ k \le k_1$ and
$$
W_1 \ge (1 + 2^{-6}) \eta. \eqno (3.5)
$$

We set
$$
l_2 = \left\lceil
\log_2  { 4 k_2  \exp\left(\left(W_1
+  2^{-l_1} \right) k_2 \right) \over  \psi }
\right\rceil,
$$
and, for $n \ge 3$, we define $l_n$ by
$$
l_{n} = \left\lceil
\log_2  { 4k_{n}  \exp\left(\left(W_{n-2}
+  2^{-l_{n-2}}\right) k_{n} \right) \over  \psi }
\right\rceil.  \eqno (3.6)
$$
Let $n \ge 2$ be an integer.
Instead of the interval $A_{k,m}$, where $k\le k_n$,
we consider the shortest dyadic interval $\hat{A}_{k,m}^n$
of the form
$$
\left( {a_1 \over 2^{l_n}},  {a_2 \over 2^{l_n}}
\right) , \,\,\,\ a_1,a_2 \in \bf{Z},
$$
which covers the interval $A_{k,m} $.
Define
$$
\hat{A}_k^n :=
\bigcup_{m=\lfloor e^{\eta k} \rfloor}^{\lceil e^{2 \eta k} \rceil}
\hat{A}_{k,m} ^n\supset A_k.
$$

Let $n \ge 1$ be an integer. We check that
$$
2 t \le l_{n+1} - l_n \le 2 h. \eqno (3.7)
$$
In particular, we have $l_{n+1} \ge l_n$, thus
$$
\hat{A}_{k,m}^n \supset\hat{A}_{k,m}^{n+1}\supset {A}_{k,m}.
$$
Here, we should note that $W_{n+1}$ is not defined yet, but
it has to satisfy
$$
\left[ W_{n+1}, W_{n+1} + {1 \over 2^{l_{n+1}}}\right] \subset J_n.
\eqno (3.8)
$$
We claim that, for any such choice of $W_{n+1}$,
we have
$$
\Leb (\hat{A}_{k,m}^{n+1}) \le 4 \Leb ({A}_{k,m})  \eqno (3.9)
$$
for every $ k \le k_{n+1}$ and for every integer $m$
such that
$$
A_{k,m} \cap J_{n-1} \neq \emptyset.  \eqno (3.10)
$$
The reason for (3.9) to be valid is as follows.
Given an integer $k$, we define $m_1 = m_1(k)$
to be the maximal $m$ for which (3.10) holds.
Then
$$
m_1 (k) \le \exp \left(\left( W_{n-1} +
2^{-l_{n-1}} \right)  k + 1 \right).
$$
From (3.6) it follows that
$$
\Leb (A_{k,m}) =  {2\psi \over mk} \ge  {2\psi \over m_1(k_{n+1})k_{n+1}} \ge  {2 \over 2^{l_{n+1}}},
$$
for $k \le k_{n+1}$, and (3.9) holds
for any possible value of $W_{n+1}$ satisfying (3.8).

 \bigskip

\noi{\it 3. Inductive assumption.}

\medskip

We describe the inductive assumption of
our version of the Peres--Schlag method.
It consists, for $n \ge 2$, of the following two points
$({\bf i_n})$ and $({\bf ii}_n)$,
that have to be satisfied by an interval $J$:

$$
\leqno ({\bf i_n})   \qquad
J \cap \hat{A}_{k}^{n-1} = \emptyset, \quad
\hbox{for every $ k \le k_{n-1}$}.
$$

$$
\leqno ({\bf ii_n}) \qquad
\Leb \bigl(J  \setminus
\bigl(\bigcup_{k\le k_n} \hat{A}_k^n\bigr)\bigr) \ge \Leb (J )/2.
$$

\medskip

Estimating
$$
\sum_{k=h+1}^{2h} \,
\sum_{m = \lfloor e^{k W_1}
\rfloor}^{\lceil  e^{k (W_1 + 2^{-l_1})} \rceil} \, {2 \psi \over mk}
\le \psi h 2^{-l_1} \le 2^{-l_1-6},
$$
our choice of $l_2$ implies that
$$
\Leb \biggl(J_1 \setminus
\biggl(\bigcup_{k\le k_2} \hat{A}_k^2 \biggr) \biggr) \ge
{\Leb (J_1) \over 2}.
$$
We have thus checked that $({\bf i_2})$ and
$({\bf ii_2})$ hold for the interval $J_1$.

\bigskip

\bigskip

\noi{\it 4. Independence shift.}

\medskip

This subsection is devoted to the proof of a key lemma
for the inductive step.

\proclaim Lemma 1.
For $n \ge 2$ and $k$ satisfying
$k_n \le k \le k_{n+1}$, we have
$$
\Leb (J_{n-1} \cap \hat{A}_k^{n+1})  < 16 \psi \Leb (J_{n-1}).
$$

\pro
It follows from (3.9) that it is enough to establish that
$$
\Leb (J_{n-1} \cap {A}_k)  < 4 \psi \Leb (J_{n-1}),   \eqno (3.11)
$$
for $k \ge k_n$. Recall that
$$
A_k =
\bigcup_{m=\lfloor e^{\eta k} \rfloor}^{\lceil e^{2 \eta k} \rceil}
\left ( {\log m \over k} -  {\psi  \over km} ,
{\log m \over k} + {\psi  \over km} \right).
$$
Let
$m_0, m_0+1, \ldots ,m_0+t = m_1$ be the integers $m$ for which
the interval $A_{k,m}$
has non-empty intersection with the segment $J_{n-1}$.
Then
$$
\Leb (J_{n-1} \cap A_k ) \le \sum_{m=m_0}^{m_1}  {2 \psi \over km}
$$
and
$$
m_0 \ge e^{W_{n-1}k} -1 \ge {e^{W_{n-1}k_n} \over 2}.
$$
We check that
$$
\max_{j=0, \ldots ,t-1}  \left(  {\log (m_0+j+1) \over k} -
 {\log (m_0+j) \over k}\right) \le
{1 \over m_0 k} \le
{2 \over k_n e^{W_{n-1} k_n}}.
$$
Furthermore, since $2^{-l_{n-1}} k_n \le 1$, we get
$$
\eqalign{
k_n e^{W_{n-1} k_n} & \ge
k_{n-1} e^{(W_{n-3} + 2^{-l_{n-3}}) k_{n-1}} e^{-2^{-l_{n-1}} k_n}
e^{W_{n-1} h}  \cr
& \ge (\psi 2^{l_{n-1}}) \cdot {1 \over 4} \cdot {64 \over \psi} \cr
& \ge 2^{l_{n-1} + 4} \ge {16 \over \Leb(J_{n-1})}, \cr}
$$
for $n \ge 4$.
We check below that the inequality
$$
k_n e^{W_{n-1} k_n} \ge {16 \over \Leb(J_{n-1})} \eqno (3.12)
$$
also holds for $n=2$ and $n=3$.

For $n = 2$, inequality (3.12) is satisfied as soon as
$$
2 h e^{2hW_1} \psi \ge 2^5 h e^{2 \eta h},
$$
that is, using (3.5), as soon as
$$
e^{2^{-5} \eta h} \psi \ge 2^4.
$$
The latter inequality is a direct consequence
of (3.1) and (3.3), provided that $t$ is sufficiently large.

For $n=3$, inequality (3.12) holds as soon as
$$
e^{3hW_2} \ge 2^5 e^{2(W_1 + 2^{-l_1}) h},
$$
which, by (3.4), holds for $t$ sufficiently large.

Consequently, for $n \ge 2$,
at least two centers of
the intervals $A_{k,m}$ are lying inside $J_{n-1}$ and
$$
\Leb (J_{n-1} )\ge
{\log(m_1/m_0) \over k} -  {2 \psi \over k_n m_0}.
$$
Thus, we get
$$
\Leb (J_{n-1} \cap A_k ) \le \sum_{m=m_0}^{m_1}  {\psi \over km}
\le \psi \Leb (J_{n-1} ) + {2\psi  \over k_n m_0}.
\eqno (3.13)
$$
Since
$$
 {1 \over k_n m_0}
\le
 {1 \over k_n( e^{W_{n-1}k_n}-1)}\le  {1 \over 2^{l_{n-1}}} =
\Leb (J_{n-1}),
$$
the lemma follows from (3.11) and (3.13). \cqfd

\bigskip

\noi{\it 5. Inductive step.}

\medskip

Let $n \ge 2$ be an integer and $J_{n-1}$
be an interval such that $({\bf i_n})$
and $({\bf ii_n})$ hold with $J = J_{n-1}$.
We consider the set
$$
J_{n-1}\setminus \left(\bigcup_{k = k_{n-1}+1}^{k_n}\hat{A}^n_k \right)
=\bigcup_{\nu = 1}^T I^\nu,
$$
where $T \ge 1$ and the $I^\nu$ are distinct intervals of the form
$$
I^\nu = \left[  {a^\nu \over 2^{l_n}},
{a^\nu+1 \over 2^{l_n}}\right].
$$
 We see that
$$
\Leb(I^\nu ) =  {1 \over 2^{l_n}}, \qquad
I^\nu \cap \hat{A}_k =\emptyset,
$$
for $\nu = 1, \ldots , T$ and for $ k\le k_n$.
For a given index $\nu$ consider the set
$$
I^\nu \setminus \bigcup_{k= k_n+1}^{k_{n+1}} \hat{A}_k^{n+1}.
$$
We see that
$$
\Leb \left( I^\nu \setminus
\bigcup_{k= k_n+1}^{k_{n+1}} \hat{A}_k^{n+1}\right)
\ge
 {1 \over 2^{l_n}} - \sum_{k=k_n+1}^{k_{n+1}}
\Leb \left( I^\nu \cap \hat{A}^{n+1}_k\right),
$$
thus
$$
\sum_{\nu = 1}^T
\Leb \left( I^\nu \setminus
\bigcup_{k= k_n+1}^{k_{n+1}} \hat{A}_k^{n+1}\right)
\ge
 {T \over 2^{l_n}} - \sum_{\nu = 1}^T  \sum_{k=k_n+1}^{k_{n+1}}
\Leb \left( I^\nu \cap \hat{A}^{n+1}_k\right).
$$
But
$$
\eqalign{
\sum_{\nu = 1}^T  \sum_{k=k_n+1}^{k_{n+1}}
\Leb \left( I^\nu \cap \hat{A}^{n+1}_k\right) & \le
\sum_{k=k_n+1}^{k_{n+1}}
\Leb \left( J_{n-1} \cap \hat{A}^{n+1}_k\right) \cr
& \le  16\psi h
\Leb \left( J_{n-1} \right), \cr}
$$
by Lemma 1.
We deduce from the inductive assumption ({\bf ii})
that
$$
\Leb(J_{n-1}) \le 2   \Leb\left( J_{n-1}\setminus
\left(\bigcup_{k = k_{n-1}+1}^{k_n}A_k \right)  \right).
$$
Furthermore, we have
$$
 {T \over 2^{l_n}} = \Leb\left( J_{n-1}\setminus
\left(\bigcup_{k = k_{n-1}+1}^{k_n}A_k \right)  \right)
$$
and, by (3.1) and for $t$ large enough,
$$
32 \psi h \le  {1 \over 2}.
$$
Consequently,
$$
\eqalign{
\sum_{\nu = 1}^T
\Leb \left( I^\nu \setminus
\bigcup_{k= k_n+1}^{k_{n+1}} \hat{A}_k^{n+1}\right)
& \ge  {1 \over 2} \cdot \Leb
\left( J_{n-1}\setminus
\left(\bigcup_{k = k_{n-1}+1}^{k_n}A_k \right)  \right) \cr
& =  {1 \over 2} \cdot
\sum_{\nu = 1}^T \Leb( I^\nu). \cr}
$$
Thus, there exists $\nu_0 = 1, \ldots , T$ such that
$$
\Leb \left( I^{\nu_0} \setminus
\bigcup_{k= k_n+1}^{k_{n+1}} \hat{A}_k^{n+1}\right)
\ge  {1 \over 2} \Leb (I^{\nu_0}).
$$
We put $J_n = I^{\nu_0}$.
We have shown that $({\bf i_n})$ and
$({\bf ii_n})$ are satisfied with $J = J_n$.

\bigskip

\noi{\it 6. Conclusion.}

\medskip

The sequence $(J_n)_{n \ge 1}$ is a decreasing
(with respect to inclusion)
sequence of non-empty compact intervals.
Consequently,
the intersection $\cap_{n \ge 1} \, J_n$ is non-empty.
By construction, if $\xi$ is in $\cap_{n \ge 1} \, J_n$,
then $\xi$ avoids every interval $\hat{A}_k^n$.
This completes the proof of Theorem 1.

We can slightly modify our construction to end up with
an uncountable intersection $\cap_{n \ge 1} \, J'_n$.
Indeed, in view of (3.7), the integer $T$
occurring in the inductive step is not too small. Thus,
at each step $n$,
we have at least two choices for the interval $J_n$,
and we let $J'_n$ be the union of two such suitable intervals.

Furthermore, a (small) positive $\delta$ being given,
we see that there are indeed at least $\lfloor (1 - \delta) T \rfloor$
suitable choices for $J_n$ at each step $n$, provided that
the value $2^{14}$ in (3.1) is replaced by a larger
number, say $\kappa(\delta)$,
depending only on $\delta$ and tending to infinity as
$\delta$ tends to $0$. Thus, we have a Cantor type
construction and the Hausdorff dimension of the resulting set
${\cal C}_{\delta}$
can be bounded from below by means of the mass distribution
principle, as was done in \cite{BuMo11}.
Replacing the value $2^{6}$ in (3.3) by $\sqrt{\kappa(\delta)}$,
it follows from a rapid calculation using (3.7) that the
Hausdorff dimension of ${\cal C}_{\delta}$ tends to $1$
as $\delta$ approaches $0$. This establishes the metrical
statement enounced after Theorem 1.
\cqfd

\vskip 5mm

\centerline{\bf 4. Proof of Theorem 3}

\vskip 6mm

We adapt the proof of Theorem 1 of \cite{Vij48}.
For simplicity, we only treat the case where $\xi = 1$.
Set $b_n = a_n + \eps$ for $n \ge 1$.
Without any loss of generality, we assume that $\eps \le 1/2$.
Let $\eta$ be a positive real number with $\eta < 1$.
Let $H$ be an integer such that $\eps H > H^{\eta} > 2$
and put $c = \lfloor H^{\eta} \rfloor - 2$.

Set $I_1 = [H + a_1, H + b_1]$. This is our Step 1.
Since
$$
(H + b_1)^2 - (H + a_1)^2 \ge 2 H^{\eta}\ge c+ 2,
$$
there is an integer $j_1$ such that $j_1, \ldots , j_1 + c$
are in the interval $[(H + a_1)^2, (H + b_1)^2]$.
For $h = j_1, \ldots , j_1 +c-1$, let $I_{2,h}$
be the interval $[\sqrt{h + a_2}, \sqrt{h + b_2}]$.
Since
$$
(H + a_1)^2 \le h + a_2 < h + b_2 \le (H + b_1)^2,
$$
the interval $I_{2,h}$ is included in $I_1$.
By construction, every real number $\xi$ in $I_{2,h}$
is such that $\{\xi\}$ and $\{\xi^2\}$ are in $[a_1, b_1]$
and $[a_2, b_2]$, respectively.
Let $E_2$ be the union of the $c$ intervals $I_{2,h}$.
This completes Step 2.

We continue this process.
Let $h = j_1, \ldots , j_1 +c-1$. Since
$$
(\sqrt{h  + b_2})^3 - (\sqrt{h + a_2})^3
\ge \bigl( (\sqrt{h + b_2})^2 - (\sqrt{h  + a_2})^2 \bigr)
\sqrt{h + a_2} \ge H \eps \ge c + 2,
$$
there is an integer $j_2$ such that $j_2, \ldots , j_2 + c$
are in the interval $[(h + a_2)^{3/2}, (h + b_2)^{3/2}]$.
For $i = j_2, \ldots , j_2 +c-1$, let $I_{2,h,i}$
be the interval $[(i + a_3)^{1/3}, (i + b_3)^{1/3}]$.
By construction, $I_{2,h,i}$ is included in $I_{2,h}$.
Proceeding in this way, we construct at Step 3 a union $E_3$ of
$c^2$ sub-intervals of $I_1$, whose elements $\xi$
have the property that $\{\xi\}, \{\xi^2\}$
and $\{\xi^3\}$ are in $[a_1, b_1]$, $[a_2, b_2]$ and $[a_3, b_3]$,
respectively.

Continuing further in the same
way, for $j \ge 4$, we construct at Step $j$ a set $E_j$
which is the union of $c^{j-1}$ closed intervals
of length approximately equal to
$$
\asymp (H^j + b_j)^{1/j} - (H^j + a_j)^{1/j}
\asymp \eps H^{-j+1} / j.
$$
Each of these intervals gives birth to $c$ intervals at the next step.
Furthermore, two different intervals at Step $j$ are separated by
at least $H^{-j+1} / j$ times an absolute positive
constant. The set
$$
{\cal C}_{\eta} = \bigcap_{j \ge 1} \, E_j
$$
is a Cantor type set, whose elements have the property that,
for $n \ge 1$,
the fractional part of their $n$-th power
lies in $[a_n, b_n]$.
The Hausdorff dimension of ${\cal C}_{\eta}$ can be bounded from
below by using the mass distribution principle,
as given, e.g., in Chapter 4 of \cite{Fal90}.
We get that the dimension of ${\cal C}_{\eta}$
is at least
equal to $(\log c)/(\log H)$ and,
since $\eta$ can be taken arbitrarily close to $1$, our theorem
is proved.  \cqfd

\vskip 8mm

\vskip 8mm

\goodbreak

\centerline{\bf References}

\vskip 5mm

\beginthebibliography{999}

\bibitem{BuMo11}
Y. Bugeaud and N. Moshchevitin,
{\it Badly approximable numbers and Littlewood-type problems},
Math. Proc. Cambridge Phil. Soc.
To appear.

\bibitem{Dub07c}
A. Dubickas,
{\it On the powers of some transcendental numbers},
Bull. Austral. Math. Soc.  76  (2007), 433--440.

\bibitem{Fal90}
K. Falconer,
Fractal Geometry : Mathematical Foundations and
Applications, John Wiley \& Sons, 1990.

\bibitem{Ko35}
J. F. Koksma,
{\it Ein mengentheoretischer Satz \"uber
die Gleichverteilung modulo Eins},
Compositio Math. 2 (1935), 250--258.

\bibitem{Mosh09a}
N. G. Moshchevitin,
{\it A version of the proof for Peres--Schlag's
theorem on lacunary sequences}.
Available at arXiv: 0708.2087v2 [math.NT] 15Aug2007.

\bibitem{Mosh09b}
N. G. Moshchevitin,
{\it Density modulo $1$ of sublacunary sequences:
application of Peres--Schlag's arguments}.
Available at arXiv:  0709.3419v2 [math.NT] 20Oct2007.

\bibitem{PeSc09}
Yu. Peres and W. Schlag,
{\it Two Erd\H os problems on lacunary sequences:
chromatic numbers and Diophantine approximations},
Bull. London Math. Soc. 42 (2010), 295--300.

\bibitem{PiSa64}
Ch. Pisot and R. Salem,
{\it Distribution modulo $1$ of
the powers of real numbers larger than $1$},
Compositio Math. 16 (1964), 164--168.

\bibitem{Vij48}
T. Vijayarhagavan,
{\it On the fractional parts of powers of a number. IV},
J. Indian Math. Soc. (N.S.)  12,  (1948), 33--39.

\endthebibliography

\vskip1cm

Yann Bugeaud \hfill Nikolay Moshchevitin

Universit\'e de Strasbourg \hfill Moscow State University

Math\'ematiques \hfill Number Theory

7, rue Ren\'e Descartes \hfill Leninskie Gory 1

67084 STRASBOURG Cedex (France) \hfill MOSCOW (Russian federation)

\medskip

{\tt bugeaud@math.unistra.fr} \hfill {\tt moshchevitin@rambler.ru}

\bye

\bye